\RequirePackage{etoolbox}
\csdef{input@path}{%
 {sty/}
 {img/}
}%
\csgdef{bibdir}{bib/}

\documentclass[ba]{imsart}
\pubyear{2016}
\volume{TBA}
\issue{TBA}
\doi{0000}
\firstpage{1}
\lastpage{2}

\usepackage{amsthm}
\usepackage{amsmath}
\usepackage{natbib}
\usepackage[colorlinks,citecolor=blue,urlcolor=blue,filecolor=blue,backref=page]{hyperref}
\usepackage{graphicx}

\startlocaldefs
\endlocaldefs
\begin{document}
\begin{frontmatter}
\title{Contributed Discussion to Bayesian Solution Uncertainty Quantification for Differential Equations by Oksana A. Chkrebtii, David A. Campbell, Ben Calderhead, and Mark A. Girolami}


\begin{aug}
\author{\fnms{William Weimin} \snm{Yoo}\thanksref{addr1}\ead[label=e1]{yooweimin0203@gmail.com}}

\address[addr1]{Mathematical Institute, Leiden University, The Netherlands \printead{e1}}


\end{aug}
\begin{abstract}
We begin by introducing the main ideas of the paper, and we give a brief description of the method proposed. Next, we discuss an alternative approach based on B-spline expansion, and lastly we make some comments on the method's convergence rate.
\end{abstract}

\begin{keyword}
\kwd{Differential equation}
\kwd{discretization uncertainty}
\kwd{B-splines}
\kwd{tensor product B-splines}
\kwd{convergence rate}
\end{keyword}
\end{frontmatter}

I would like to congratulate the authors for such an interesting research. The Bayesian method with the probabilistic solver introduced is highly innovative and practical. The various examples presented in the paper show the wide applicability of the proposed method. However, I do find the title a bit of a misnomer, as I initially thought that the authors are constructing credible sets for the fixed but unknown solution $u^{*}$ of the differential equation.

The inverse problem that the authors are trying to solve, in its most basic form is this: Suppose you have observations $Y=Au+\varepsilon$, where $\varepsilon$ is some normal errors and $u$ follows $u_t=f(t,u,\theta)$. Here, $A$ is a known transformation from the state space $u$ to the observation space $Y$, $u_t$ is the first order derivative with respect to its argument $t$, $f$ is the known form of the differential equation, and $\theta$'s are the equation's parameters. The method proposed consists of two steps, with one nested within the other. First, solve for $u$ probabilistically to obtain a discretized solution at some grid points. Then we embed these discretized version of $u$ in a Bayesian hierarchical framework to estimate $\theta$. To model discretization uncertainty associated with using only $u$ evaluated at grid points, the authors endow priors based on Gaussian process jointly on $u$ and $u_t$, where the covariance function is constructed by convolving kernels.

There is an alternative and perhaps a conceptually easier way to achieve the same result. We can first represent $u$ by a B-spline series, i.e., $u(t)=\sum_{j=1}^J\vartheta_jB_{j,q}(t)$ with $B_{j,q}(\cdot)$ denoting the $j$th B-spline of order $q$, and we endow the coefficients $\vartheta_j$'s with normal priors. Here, the number of basis $J$ plays the role of $1/\lambda$, where $\lambda$ is the length-scale parameter defined in the paper. It turns out that the first derivative of this $u$ is another B-spline series $u_t(t)=\sum_{j=1}^{J-1}\vartheta_j^{(1)}B_{j,q-1}(t)$ where $\vartheta_j^{(1)}$ is some weighted first order finite difference of the $\vartheta_j$'s ((4.23) of \citet{lschumaker}). Therefore, $u$ and $u_t$ are jointly normal and their associated covariance matrices are banded due the support separation property of B-splines. To enforce the given initial condition, we can condition the joint prior $(u,u_t)$ on $u^{*}(0)$.

Moreover, this approach can be generalized to the partial differential equation case, where we take tensor product of B-splines to model both the spatial and temporal components, i.e., $u(x,t)=\sum_{j_1=1}^J\sum_{j_2=1}^J\theta_{j_1,j_2}B_{j_1,q}(x)B_{j_2,q}(t)$. As in the univariate case, partial derivatives of tensor product B-splines will be another tensor-product B-splines ((3.2) of \cite{yoo2016}). Hence we will obtain the same Gaussian process prior for $u$ and all its mixed partial derivatives if we endow normal priors on the coefficients. As before, we enjoy some simplification in computing the covariance matrices because they are banded.

In addition, I would like to comment on the effect of grid point distribution on the convergence rate of the proposed algorithm. Intuitively, one would expect that the grid points should be chosen roughly uniformly across the domain $[0,L]$. Suppose we choose grid points $\{t_1,t_2,\dotsc,t_N\}$ and we further assume that they are quasi-uniform, i.e., $h/\min_i(t_i-t_{i-1})\leq C$ for some constant $C>0$ with $h=\max_i(t_i-t_{i-1})$. In other words, the max grid increment is of the same order as the min grid increment. Then it follows that $h$ is of the order of $1/N$ and by Theorem 1 of \citet{discuss}, the rate of convergence is $O(N^{-1})$. Therefore for quasi-uniform grids (which includes uniform discrete grids), increasing the number of grid points will result in more accurate solution.

The paper under discussion \citet{discuss} makes significant contribution to the new field of probabilistic numerics. I have learnt a great deal by reading this paper, and I hope that there will be more papers in uncertainty quantification for differential equation models in the future.

\bibliographystyle{ba}
\bibliography{discussref}
\end{document}